	\documentclass[12pt, a4paper]{article}
	
 	\usepackage{latexsym}
	\usepackage{amsmath}	
	\usepackage{epsfig}	
	\usepackage{times}	
	\usepackage{amssymb}

 	\newtheorem{thr}{Theorem}[section]
 	\newtheorem{prop}{Proposition}[section]
 	
	\newtheorem{cor}{Corollary}[section]
	\newtheorem{lem}{Lemma}[section]
 	\newcommand{\F}{\mathbb{F}}

\begin{document}
 %
 	\centerline{\Large{\bf On the Divisibility of Trinomials by Maximum}}
 	\centerline{}
 	\centerline{\Large{\bf  Weight Polynomials over $\F_2$}}
 	\centerline{}
 	\centerline{}
 	\centerline{\textsuperscript{a}Ryul Kim, Ok-Hyon Song and Myong-Hui Ri}
 	\centerline{}
 	{\small \centerline{ Faculty of Mathematics,  \textbf{Kim Il Sung} University, 
	Pyongyang, D.P.R Korea}}
	{\small \centerline{{e-Mail address : \textsuperscript{a}ryul\_kim@yahoo.com}}}
	\centerline{}
	\centerline{}
%
 %
	\begin{abstract}
	Divisibility of trinomials by given polynomials over finite fields has been studied and used to construct orthogonal arrays
	in recent literature.
	Dewar et al.\ (Des.\  Codes Cryptogr.\ 45:1-17, 2007) studied the division of
	trinomials by a given pentanomial over $\F_2$ to obtain the orthogonal arrays of strength at least 3, and finalized their paper 
	with some open questions. One of these questions is concerned with generalizations to the polynomials with
	more than five terms.
	In this paper, we consider the divisibility of trinomials by a given maximum weight polynomial over $\F_2$ and 
	apply the result to the construction of the orthogonal arrays of strength at least 3.
	\end{abstract}
%
%
%
%
	\section{Introduction}
	Sparse irreducible polynomials such as trinomials over $\F_2$ are widely used to perform 
	arithmetic in extension fields of $\F_2$ due to fast modular reduction. In particular, primitive trinomials
	and maximum-length shift register sequences generated by them play an important role in various applications such as 
	stream ciphers (see \cite{go1,jam}). But even irreducible trinomials do not exist for every degree. 
	When a primitive (respectively irreducible) trinomial of a given degree does not exist, an almost 
	primitive (respectively irreducible) trinomial, which is a reducible trinomial with primitive
	(respectively irreducible) factor, may be used as an alternative \cite{bre}.
	 This encouraged the researchers to study divisibility of trinomials by primitive 
	or irreducible polynomials \cite{che,go2,kim}. The divisibility of trinomials by primitive polynomials is also 
	related to orthogonal arrays.

	Let $f$ be a polynomial of degree $m$ over $\F_2$ and let $a=(a_0,a_1,\cdots)$ be a
	shift-register sequence with characteristic polynomial $f$. Denote by $C_n^f$ the
	set of all subintervals of this sequence with length $n$, where $m<n \leq 2m$,
	together with the zero vector of length $n$. 	
	 Munemasa \cite{mun} observed that very few trinomials of degree at most $2m$ are divisible by 
	a given primitive trinomial of degree $m$ and proved that if $f$ is a primitive trinomial satisfying
	certain properties, then $C_n^f$ is an orthogonal array of strength 2 having the property
	of being very close to an orthogonal array of strength 3.
	Munemasa's work was extended in \cite{dew}. The authors considered the divisibility of a trinomial of degree at most $2m$ 
	by a given pentanomial $f$ of degree $m$ and obtained the orthogonal arrays of strength 3. 
	They suggested some open questions in the end of their paper. One of them is to extend the results to 
	finite fields other than $\F_2$. In this regard, Panario et al.\ \cite{pan} characterized the divisibility of binomials and 
	trinomials over $\F_3$. Another question in \cite{dew} is related to extend the results to the polynomials with
	more than five terms. In this paper we analyze the division of trinomials by a maximum weight polynomial over $\F_2$.
	
	In the theory of shift register sequences it is well known that the lower the weight, i.e. the number 
	of nonzero coefficients of the characteristic 
	polynomial of shift register sequence, is, the faster is the generation of the sequence.
	But Ahmadi et al.\ \cite{ahm} point out the advantage of maximum weight polynomials over $\F_2$ in the implementation 
	of fast arithmetic in extension fields. We show that no trinomial of degree at most $2m$ is divisible by a given 
	maximum weight polynomial $f$ of degree $m$, provided that $m>7$. Using this result we can also 
	obtain the orthogonal arrays of strength at least 3. The rest of the paper is organized as follows. 
	In Section 2, some basic definitions and results are given and in Section 3, some properties of maximum weight polynomials 
	and shift register sequences generated by them are mentioned. We focus on the divisibility of trinomials by maximum 
	weight polynomials in Section 4, and conclude in Section 5.		
%
%
 %
%
	\section{Preliminaries}

	A \textit{period} of a nonzero polynomial $f(x) \in \F_q[x]$ with $f(0) \neq 0$  is the least positive integer 
	$e$ for which $f(x)$ divides $x^e-1$. A polynomial $f(x)\in \F_q[x]$ is called \textit{reducible} if it has 
	nontrivial factors; otherwise \textit{irreducible}. A polynomial $f(x)$ of degree $m$ is called 
	\textit{primitive} if it is irreducible and has period $2^m-1$. The \textit{reciprocal polynomial} of 
	$f(x)=a_mx^m+a_{m-1}x^{m-1}+\cdots+a_1x+a_0 \in \F_q[x]$ with $a_m \neq 0$ is defined by
	\begin{equation*}
	f^*(x)=x^mf(1/x)=a_0x^m+a_1x^{m-1}+\cdots+a_{m-1}x+a_m.
	\end{equation*}
	We refer to \cite{lid} for more information on the polynomials over finite fields. Throughout this 
	paper we only consider a binary field $\F_2$ and all the polynomials are assumed to be in $\F_2[x]$, 
	unless otherwise specified. 

	A \textit{shift-register sequence} with
	characteristic polynomial $f(x)=x^m+\sum_{i=0}^{m-1}c_ix^i$ is the sequence
	$a=(a_0,a_1,\cdots)$ defined by the recurrence relation
	\begin{equation*}
	a_{n+m} = \sum_{i=0}^{m-1}c_i a_{i+n}
	\end{equation*}
	\noindent for $n \geq 0$.
	
	A subset $C$ of $\F_2^n$ is called an \textit{orthogonal array} of
	strength $t$ if for any $t-$ subset $T=\{i_1,i_2,\cdots,i_t\}$ of $\{1,2,\cdots,n\}$
	and any $t-$tuple $(b_1,b_2,\cdots,b_t) \in \F_2^t$, there exists exactly
	$|C|/2^t$ elements $c=(c_1,c_2,\cdots,c_n)$ of $C$ such that $c_{i_j}=b_j$
	for all $1 \leq j \leq t$\cite{mun}. From the definition, if $C$ is an orthogonal array
	of strength $t$, then it is also an orthogonal array of strength $s$ for all
	$1 \leq s \leq t$.
%
%
	
	The next theorem, due to Delsarte, relates orthogonal arrays to linear codes.
	\begin{thr} \textnormal{(\cite{del})}
	Let $C$ be a linear code over $\F_q$. Then $C$ is an orthogonal array of
	maximum strength $t$ if and only if $C^{\bot}$, its dual code, has minimum
	weight $t+1$.
	\end{thr}

	Munemasa \cite{mun} described the dual code of the code generated by
	a shift-register sequence in terms of multiples of its primitive characteristic
	polynomial and Panario et al.\ \cite{pan} generalized this result as follows by removing 
	the primitiveness condition for the characteristic polynomial. 
%
%
	\begin{thr} \textnormal{(\cite{pan})}
	Let $a=(a_0,a_1,\cdots)$ be a shift register sequence with minimal polynomial $f$, 
	and suppose that $f$ has degree $m$ with $m$ distinct roots. Let $\rho$ be the period of $f$ 
	and $2 \leq n \leq \rho$. Let $C_n^f$ be the set of all subintervals of the shift register 
	sequence $a$ with length $n$, together with the zero vector of length $n$. 
	Then the dual code of $C_n^f$ is given by
	\begin{equation*}
	(C_n^f)^{\bot} = \{(b_1,\cdots,b_n) : \sum_{i=0}^{n-1}b_{i+1}x^i
	\textnormal{ is divisible by } f \}.
	\end{equation*}\	
	\end{thr}	

	A \textit{maximum weight polynomial} is a degree-$m$ polynomial of weight $m$ (where $m$ is odd) over $\F_2\cite{ahm}$, namely,
	 \begin{equation*}
	f(x)=x^m+x^{m-1}+\cdots+x^{l+1}+x^{l-1}+\cdots+x+1
	=\frac{x^{m+1}+1}{x+1}+x^l.
	\end{equation*}
	If you take
	 \begin{equation*}
	g(x)=(x+1)f(x)=x^{m+1}+x^{l+1}+x^l+1,
	\end{equation*}
	then the weight of $g(x)$ is 4, and its middle terms are consecutive, so reduction using $g(x)$ instead of $f(x)$
	is possible and can be effective in the arithmetic of an extension field $\F_{2^m}$ as if the reduction 
	polynomial were a trinomial or a pentanomial.
	This fact motivated us to consider the divisibility of trinomials by maximum weight polynomials.
%
%
%
%
	\section{Character of shift register sequence generated by a maximum weight polynomial}
	In this section we state a simple property of maximum weight polynomials and characterize the shift register 
	sequences generated by them. 
%
%
	\begin{prop}
	Let $f(x)=x^m+x^{m-1}+\cdots+x^{l+1}+x^{l-1}+\cdots+1 \in \F_2[x]$. If $f(x)$ is irreducible, 
	then $\textnormal{gcd}(m,l)=1$.	
	\end{prop}
	\textit{Proof.} Suppose $\textnormal{gcd}(m,l)=d>1, m=m_1d$ and $ l=l_1d$. Then we have
	\begin{eqnarray*}
	g(x) & := & (x+1)f(x) = x^{m+1}+x^{l+1}+x^l+1\\
			& = & x^{l+1}(x^{m-l}+1)+(x^l+1) =x^{l+1}(x^{m_1d-l_1d}+1)+(x^{l_1d}+1)\\
			& = & x^{l+1}(x^{d(m_1-l_1)}+1)+(x^{l_1d}+1).
	\end{eqnarray*}
	So $(x^d+1)/(x+1)$ is a factor of $f(x)$, which means $f(x)$ is reducible. $\Box$
%
%
	\begin{prop}
	Let $f(x)=x^m+x^{m-1}+\cdots+x^{l+1}+x^{l-1}+\cdots+1 \in \F_2[x]$ be a primitive polynomial and
	\begin{equation*}
	a_{n+m}=\sum_{i=0}^{m-1}a_{n+i} +a_{n+l} (n \geq 0)
	\end{equation*}
	be a shift-register sequence with characteristic polynomial $f$. Then for all positive integer $n$,
	\begin{eqnarray*}
	a_{n+m}=a_{n-1}+a_{n-1+l}+a_{n+l}.
	\end{eqnarray*}
	\end{prop}
	\textit{Proof.} Since $f(x)$ is the characteristic polynomial of $(a_0,a_1,\cdots)$,
	we get $a_l=a_0+a_1+\cdots+a_m$ where $a_0,a_1,\cdots,a_{m-1}$ are initial values not
	all of which are zero. We use induction on $n$.\\
	\indent If $n=1$,
	\begin{eqnarray*}
	a_{m+1} & = & a_1+\cdots+a_l+a_{l+2}+\cdots+a_m \\
		& = &a_0+ (a_0+\cdots+a_l+a_{l+1}+a_{l+2}+\cdots+a_m)+a_{l+1} \\
		& = & a_0+a_l+a_{l+1}.
	\end{eqnarray*}
	Now assume that the equation $a_{n+m}=a_{n-1}+a_{n-1+l}+a_{n+l}$ holds true
	for all positive integers less or equal to $n$. Then,
	\begin{eqnarray*}
	a_{m+n+1} & = & a_{n+1}+\cdots+a_{n+l}+a_{n+l+2}+\cdots+a_{n+m} \\
		& = & (a_0+\cdots+a_m)+(a_0+\cdots+a_n)+a_{n+l+1} \\
		& & +(a_{m+1}+\cdots+a_{m+n}) \\
		& = & a_l+(a_0+\cdots+a_n)+a_{n+l+1} +(a_0+a_l+a_{l+1})\\
		& & +(a_1+a_{l+1}+a_{l+2})+\cdots+(a_{n-1}+a_{l+n-1}+a_{l+n})\\
		&=&a_n+a_{l+n}+a_{n+l+1}
	\end{eqnarray*}
	This completes the proof. $\Box$
%
%
 %
%
	\section{Divisibility of trinomials by maximum weight polynomials}
	In this section we consider the divisibility of trinomials by
	maximum weight polynomials, provided that the degree of the trinomial does not
	exceed double the degree of the maximum weight polynomial. 
	Let $f(x)=x^m+x^{m-1}+\cdots+x^{l+1}+x^{l-1}+\cdots+1 \in \F_2[x]$ and suppose that $f(x)$ divides 
	a trinomial $g(x)$ with
	\begin{equation*}
	g(x)=f(x)h(x)=(x^m+x^{m-1}+\cdots+x^{l+1}+x^{l-1}+\cdots+1)\cdot\sum_{k=0}^{t}x^{i_k},
	\end{equation*}
	where $x^{i_k}$s are the non-zero terms of $h(x)$ and $0=i_0<i_1<\cdots<i_t$.
	The above equation can be illustrated as in Fig. 1.
%
%

	\vspace{0.5cm}	

	\newcounter{cms}
	\setlength{\unitlength}{0.35mm}
	
	\begin{picture}(10,15)
	\small
	\put(20,0){$m$}
	\put(40,0){$m-1$}
	\put(80,0){$\cdots$}
	\put(110,0){$l+1$}
	\put(140,0){$(l)$}
	\put(160,0){$l-1$}
	\put(190,0){$\cdots$}
	\put(220,0){0}
	\put(350,0){$i_t$}	

	\put(50,-15){$m$}
	\put(70,-15){$m-1$}
	\put(110,-15){$\cdots$}
	\put(140,-15){$l+1$}
	\put(170,-15){$(l)$}
	\put(190,-15){$l-1$}
	\put(220,-15){$\cdots$}
	\put(250,-15){0}
	\put(350,-15){$i_{t-1}$}	

	\put(80,-40){$\ddots$}
	\put(160,-40){$\ddots$}
	\put(240,-40){$\ddots$}
	
	\put(80,-60){$m$}
	\put(100,-60){$m-1$}
	\put(140,-60){$\cdots$}
	\put(170,-60){$l+1$}
	\put(200,-60){$(l)$}
	\put(220,-60){$l-1$}
	\put(250,-60){$\cdots$}
	\put(280,-60){0}
	\put(350,-60){$i_1$}	

	\put(110,-75){$m$}
	\put(130,-75){$m-1$}
	\put(170,-75){$\cdots$}
	\put(200,-75){$l+1$}
	\put(230,-75){$(l)$}
	\put(250,-75){$l-1$}
	\put(280,-75){$\cdots$}
	\put(310,-75){0}
	\put(350,-75){$i_0$}	

	\put(10,-80){\line(1, 0){350}}
	\put(15,-75){+}
	
	\put(20,-100){\framebox(15,15){}}
	\put(140,-100){\framebox(15,15){}}
	\put(305,-100){\framebox(15,15){}}	

	\end{picture}
		
	\vspace{3.7cm}
	\small
	\centerline{{\bf Fig. 1} An illustration of equation $g(x)=f(x) \sum_{k=0}^{t}x^{i_k}$}
	
	\vspace{0.5cm}
	\normalsize

	\noindent Here $(l)$ stands for the missing terms. We adopt the same terminology as in \cite{dew,pan}. 
	In particular, if the sum of coefficients in the same column of Fig.\ 1 is 0, 
	then we write that the corresponding terms $x^i$ \textit{cancel} and if the sum is 1 then we say that one of the 
	corresponding terms is \textit{left-over}. 
	The proof of our main results will be done with Fig.\ 1. Since the most top-left term $m+i_t$ and the most bottom-right term 
	$0+i_0$ are trivial left-over terms, we have only one left-over term undetermined. Below a \textit{left-over term} 
	means the left-over term which is neither $m+i_t$ nor $0+i_0$. And we always assume that 
	$m+i_0$ is in the same column as $s+i_t, 0  \leq s  \leq m-1$ and denote the number of terms in $h(x)$ as $N$.

%
%
	
	\begin{lem}
	Let  $f(x)=x^m+x^{m-1}+\cdots+x^{l+1}+x^{l-1}+\cdots+1 \in \F_2[x]$ and $g(x)$ be a trinomial 
	of degree at most $2m$ divisible by $f(x)$ with $g(x)=f(x)h(x)$. Then $N$ equals to 3 or 5.
	\end{lem}

	\noindent \textit{Proof.} Since $g(x)$ is a trinomial and $f(x)$ has an odd number of terms, 
	$h(x)$ also has an odd number of terms, that is, $t$ is even. Suppose that $N$ is greater or equal to 7. 
	If $s \geq l$ then for every even number $k$, $m+i_{t-k}$ is a left-over term. 
	Since $t \geq 6$, we have more than 2 left-over terms which contradicts the assumption. \\
	\indent Consider the case of $s<l$. First assume that there exists a unique left-over term to the left of $m+i_0$. 
	It is sufficient to show $l \geq 3$ 	because if so, $0+i_2$ is an extra left-term which leads to a contradiction. 
	Observe a position $l+i_t$. 
	If $l+i_t \geq m+i_{t-2}$ then clearly $l \geq i_{t-2}-i_0 \geq 4$, so we have done. Assume that $l+i_t<m+i_{t-2}$. 
	Then $l+i_t \geq m+i_{t-4}$ because if not, then $m+i_{t-2}$ and $m+i_{t-4}$ are left-over terms. 
	Thus we have $l \geq i_{t-4}-i_0$. If $l+i_t>m+i_{t-4}$ then $l>2$ and if $l+i_t=m+i_{t-4}$ then 
	 $i_{t-4}-i_0>2$ because if $i_{t-4}-i_0=2$ then $m+i_{t-5}=l+i_{t-1}$ and so an extra left-over term appears. \\
	\indent Next assume that there is no left-over term to the left of $m+i_0$. Then it is clear that $m+i_{t-2}=l+i_t$ and 
	$l \geq i_{t-2}-i_0 \geq 5$ hence $0+i_2$ and $0+i_4$ are left-over terms; contradiction. $\Box$
	
%
%
	\begin{lem}
	Under the same condition as in Lemma 1, if $s<l$ then $m+i_0$ cannot be a left-over term.
	\end{lem}
	
	\noindent \textit{Proof.} Assume that $m+i_0$ is a left-over term. Then all the 
	remaining terms in other columns must cancel and by Lemma 1 $N=3$ or $N=5$. 
	If  $N=3$, then $l+i_1>m+i_0$ from $s<l$ and thus an 
	extra left-over term occurs in the column of $l+i_1$. Now assume that $N$ is 5. 
	We see easily $l+i_t=m+i_{t-2}$ and $i_t-i_{t-1}=1$. If there is an extra left-over term to 
	the left of $m+i_0$, then we have done. If there is no any extra left-over term to the left of $m+i_0$, then $i_2-i_1=2$ 
	because if $i_2-i_1=1$ then $m+i_1=l+i_{t-1}$ and so $m+i_1$ is an extra left-over term and 
	if $i_2-i_1>2$ then $l-2+i_t=l-1+i_{t-1}=m-2+i_2$ and so $l-2+i_t$ is an extra left-over term. 
	Then from the condition $i_t \leq m$, 
	it follows $l \geq 3$ and thus $0+i_2$ is an extra left-over term; contradiction. $\Box$

%
%

	\begin{thr}
	Let  $f(x)=x^m+x^{m-1}+\cdots+x^{l+1}+x^{l-1}+\cdots+1 \in \F_2[x]$.
	If $g(x)$ is a trinomial of degree at most $2m$ divisible by $f(x)$ with $g(x)=f(x)h(x)$, then \\
	\indent \textnormal{1)} $f(x)$ is one of the polynomial exceptions given in Table 1.	\\
	\indent \textnormal{2)} $f(x)$ is the reciprocal of one of the polynomials listed in the previous item.
	\end{thr}

%
%

	\small
	\centerline{
	\textbf{Table 1.} Table of polynomial exceptions
	}
	\vspace{0.3cm}
	\centerline{
	\begin{tabular}{l l l l}
	\hline
	No& $g(x)$ & $f(x)$ & $h(x)$ \\	
	\hline
	1 & $x^5+x^4+1$ & $x^3+x+1$ & $x^2+x+1$ \\
	2 & $x^6+x^4+1$ & $x^3+x^2+1$ & $x^3+x^2+1$ \\
	3 & $x^9+x^7+1$ & $x^5+x^3+x^2+x+1$ & $x^4+x+1$ \\
	4 & $x^7+x^5+1$ & $x^5+x^4+x^3+x+1$ & $x^2+x+1$ \\
	5 & $x^8+x^5+1$ & $x^5+x^4+x^3+x^2+1$ & $x^3+x^2+1$ \\
	6 & $x^{14}+x^{13}+1$ & $x^7+x^6+x^5+x^4+x^3+x+1$ & $x^7+x^5+x^2+x+1$ \\
	7 & $x^{13}+x^{10}+1$ & $x^7+x^6+x^5+x^4+x^3+x^2+1$ & $x^6+x^5+x^3+x^2+1$ \\
	\hline
	\end{tabular}
	}
	\normalsize
	\vspace{0.5cm}

	\noindent \textit{Proof.} We divide into three cases: $s>l$ or $s=l$ or $s<l$. \\
	\noindent {\bf Case 1 : } $s>l$. \\
	\indent Since $h(x)$ has an odd number of terms, $s \leq m-2$ and $m+i_0$ is a left-over term, 
	hence all the remaining terms in other columns must cancel. There is no missing term to the left of $s+i_t$, 
	and therefore $m+i_{t-2}$ is a left-over term. This means $i_0=i_{t-2}$, namely, $N=3$.
	Since $m-1+i_0$ must cancel, $s=l+1$ and  $m-2+i_0$ cancels up automatically from $i_t-i_{t-1}=1$. 
	We see easily that $l=1$ or $m-3+i_0$ is a missing term because $m-3+i_0$ must cancel up.
	If $l=1$, then clearly $m=5$ and we get the 5th polynomial in Table 1. 
	If  $m-3+i_0$ is a missing term, then $l=m-3$. Since $l-1+i_0$ must cancel up, $l$ must equal to 2 and so we get 
	the 4th polynomial in Table 1.\\
	\noindent {\bf Case 2 : } $s=l$. \\
	\indent In this case, $m+i_0$ cannot be a left-over term because the number of non-zero terms in column of $m+i_0$ is even. 
	If there is a unique left-over term to the left of $m+i_0$, then it must be $m-1+i_t$ or $m+i_2$. \\
	\noindent {\bf Case 2.1 : } $m-1+i_t$ is a unique left-over term to the left of $m+i_0$. \\
	\indent Clearly $i_{t-1}=i_t-2$. If $N=3$ then $m-1+i_0$ is an extra left-term 
	and if $N=5$ then $m+i_{t-2}$ is so. This contradicts to the assumption. \\
	\noindent {\bf Case 2.2 : } $m+i_2$ is a unique left-over term to the left of $m+i_0$. \\
	\indent This is the case of $N=5$ and $i_t-i_{t-1}= i_2-i_1=1$.
	$m-1+i_0$ cancels automatically because $m-1+i_0=l+i_{t-1}$. Thus we have only two possible cases: 
	$l=1$ or $l \neq 1, l+i_2=m-2+i_0$. Assume that $l=1$ then $m-3+i_0$ must be in the column of $l+i_2$ 
	and $m-5+i_0$ must cancel with $0+i_1$ so we get the 7th polynomial in Table 1. And assume 
	that $l \neq 1, l+i_2=m-2+i_0$ then $i_{t-1}-i_2=1$ and observing $m-4+i_0$ implies 
	that $m-4=l, l-3 \neq 0$ or $m-4>l, l=3$. In these two cases we have an extra left-over term $l-2+i_0$; contradiction.\\
	\noindent {\bf Case 2.3 : } There is no left-over term to the left of $m+i_0$. \\
	\indent It is obvious that $N=3$ and $i_t-i_1=1$. If $i_1-i_0>3$ then 
	we have two left-over terms among $j+i_0 (1 \leq j \leq 3$). Hence $i_1-i_0$ is less or equals to 3. 
	Examining all cases for $i_1-i_0$ we get the reciprocals of the 1st, 3rd and 4th polynomials in Table 1. \\
	\noindent {\bf Case 3 : } $s<l$. \\
	\indent By lemma 2, $m+i_0$ is not a left-over term. So there exists $z (1 \leq z \leq t-1)$  
	such that $m+i_0=l+i_z$. \\
	\noindent {\bf Case 3.1 : } $m+i_0=l+i_{t-1}$. \\
	\indent Clearly we have $i_{t-1} \geq i_t-3$. First assume that $i_{t-1}=i_t-3$. 
	Then $l$ equals to $m-1$ or $m-2$.
	If $l=m-1$, then $l-1+i_t=m-2+i_t$ is a left-over term so $l-3+i_t=l+i_{t-1}=m+i_0$ and $h(x)$ has three terms.
	Since the unique left-over term has already been determined, $0+i_t=l-1+i_{t-1}=l+i_0$ and 
	we get the 3rd polynomial in Table 1. 
	If $l=m-2$, then $m-1+i_t$ is a left-over term and $m+i_0$ must cancel with $0+i_t$ 
	which means $i_1-i_0=2$ and $l=3$. But then $1+i_0$ appears as an extra left-over term; contradiction.\\
	\indent Next assume that $i_{t-1}=i_t-2$. When $l \neq m-1$, $m-1+i_t$ is a left-over term and $l \leq m-3$ 
	because if $l=m-2$ then $m+i_{t-1}$ is an extra left-over term. $l+i_t$ must cancel with $m+i_{t-2}$ 
	and in fact $N$ is 5. Thus $i_2-i_1=1$. 
	By the condition $m+i_0=l+i_{t-1}$, we have $i_1-i_0=1$. Since $m-1+i_0$ must cancel up, $l-2=0$ or $m-3=l$.
	If $l-2=0$ then we get the 6th polynomial in Table 1 and the equation $m-3=l$ leads to a contradiction due to 
	an extra left-over term in column of $l-3+i_0$. When $l=m-1$, clearly $N$ is 3 from the conditin $l+i_{t-1}=m+i_0$.
  	By research of possible values of $l$ we get the reciprocals of the 2nd and 5th polynomials in Table 1.\\
	\indent Next assume that $i_{t-1}=i_t-1$. If $N=5$ then $m+i_{t-2}$ is a left-over term and $i_{t-2}-i_1=1$, 
	hence an extra left-over term occurs in the comumn of $l+i_t$. Thus $N$ is 3. Since $l-1+i_t=l+i_{t-1}=m+i_0$, 
	$l+1+i_{t-1}$ is a left-over term. If $m-1 \neq l$, then $l-1=0$ from consideration of $m-1+i_0$ 
	and therefore we get the 2nd polynomial in Table 1. If $m-1=l$, then $l-1$ cannot be zero, 
	so we get the 1st polynomial in Table 1.\\
	\noindent {\bf Case 3.2 : } $m+i_0=l+i_2$. \\
	\indent In this case $N$ is 5 and clearly $2 \leq l \leq m-2$. Observe a column of $l+i_t$.\\
	\noindent {\bf Case 3.2.1 : } $m+i_2<l+i_t$. \\
	\indent We have a left-over term in the column of $l+i_t$ and $i_t-i_{t-1}=1$. Then $m+i_2$ must cancel with $l-1+i_t$ 
	and also $i_2-i_1-1$. By the condition $l+i_2=m+i_0$, $m-1+i_0$ must cancel with $l+i_1$. 
	From $i_t \leq m$ we have $l \geq 3$ and $i_1-i_0=1$ because if not, then $1+i_0$ is an extra left-over term. 
	Hence $l$ equals to $m-2$. Since $m-1+i_0$ must cancel up, $l-4 \neq 0$. Observing the term $l-1+i_0$, we see that 
	$l-5=0$ and then $l-2+i_0$ appears as an extra left-over term; contradiction.\\
	\noindent {\bf Case 3.2.2 : } $m+i_2=l+i_t$. \\
	\indent Assume that $m-1+i_t$ is a left-over term. Then clearly $l<m-2$ and $i_t-i_{t-1}=2$. If $i_2-i_0=2$, then 
	$m+i_0$ must concel with $l+i_{t-1}$ which contradicts to the condition $m+i_0=l+t_{t-2}$. And if $i_2-i_0>2$, then 
	an extra left-over term occurs in the column of $l+1+i_t$ or $l+2+i_t$ which again leads to a contradiction.\\
	\indent Now assume that $m-1+i_t$ is not a left-over term. Then $i_t-i_{t-1}=1$ and $m+i_1$ cancels with $l+i_{t-1}$ 
	or $m+i_1<l+i_{t-1}$. If $m+i_1$ cancels with $l+i_{t-1}$ then $m+i_1$ is a left-over term and $i_2-i_1=1$. From 
	$i_t \leq m$, we have $0 \leq l-2$. Since if $i_1-i_0 \geq 2$ then $1+i_0$ is an extra left-over term, 
	$i_1-i_0=1$ and $l=m-2=4$. Then $l+2+i_0$ appears as an extra left-over term; contradiction. 
	If $m+i_1<l+i_{t-1}$ then $m+i_1$ must cancel with $m-2+i_2$ or $m-3+i_2$. Briefly considering as above, 
	we arrive at a contradiction in both cases. \\
	\noindent {\bf Case 3.2.3 : } $m+i_2>l+i_t$. \\
	\indent You shall see that $l \leq m-3, i_t-i_{t-1}=1$ and $m+i_2$ is a left-over term. 
	Since $m-1+i_2$ must cancel, $m-1+i_2=l+i_t$ or $m-1+i_2=m+i_1$. In the first case $i_2-i_1=3$ 
	because $l+i_{t-1}=m-2+i_2=l-1+i_t$. Since $m+i_1<l+i_{t-1}$, $l$ is greater or equals to 3. 
	If $i_1-i_0>1$ then $1+i_0$ is an extra left-over term and if $i_1-i_0=1$ then $l=3$ and $m-2+i_0$ is an extra 
	left-over term, which leads to a contradiction. In the second case we have $l+i_t=m+i_0$; contradiction. \\
	\noindent {\bf Case 3.3 : }  $m+i_0=l+i_1$. \\
	\indent In this case we have $l \geq 3$ from $i_t \leq m$. First assume that $1+i_0$ is a left-over term. 
	Then clearly $i_1-i_0=2, l+i_0=0+i_2$ and $l+1+i_0=l-1+i_1=1+i_2=0+i_{t-1}$. 
	Since $l+2+i_0=l+i_1=2+i_2=1+i_{t-1}=0+i_t$, we have $m=l+2$. Then from $5+i_2=4+i_{t-1}=3+i_t$, we have 
	$l=5$ which corresponds the reciprocal of the 6th polynomial in Table 1. \\
	\indent Next assume that $1+i_0$ is not a left-over term. Then $i_1-i_0=1, l=m-1$ and $0+i_2$ is a left-over term 
	because if not, then $0+i_2=l+i_0$ and thus $N=3$ which is the case mentioned above. 
	Considering the first and last terms in every rows, we have the following equations:\\
	\centerline{$i_{t-1}-i_2=1, 0+i_t=l+i_0, l+i_2>m+i_1, i_2-i_1=2,$}
	\centerline{$0+i_t=l+i_0, i_t-i_{t-1}=2.$}
	\noindent This implies the reciprocal of the 7th polynomial in Table 1. $\Box$

	Note that every polynomial $f(x)$ listed in Table 1 has degree less than 8. From this fact we can immediately 
	get the following corollary.
	\begin{cor}
	Let $f(x)$ be a maximum weight polynomial 
	of odd degree $m$ greater than 7 and g(x) be a trinomial of degree at most $2m$. Then $g(x)$ 
	is not divisible by $f(x)$.
	\end{cor}

	Combining these facts with Theorem 1 and Theorem 2, we get the
	following corollary on orthogonal arrays of strength 3.
	\begin{cor}
	Let $f(x)$ be a primitive maximum weight polynomial 
	of odd degree $m$ greater than 7. If $m \leq n \leq 2m$, then
	$C_n^f$ is an orthogonal array of strength at least 3.
	\end{cor}
	\section{Conclusion}
	In this paper, we analyzed the divisibility of trinomials by
	maximum-weight polynomials over $\F_2$ and used the result to obtain the orthogonal
	arrays of strength 3. More precisely, we showed that if $f(x)$ is a
	maximum-weight polynomial of degree $m$ greater than 7, then
	$f(x)$ does not divide any trinomial of degree at most $2m$.
	Our work gives a partial answer to one of the questions posted in
	\cite{dew}. As anticipated in \cite{dew,pan}, one seems to
	need some new techniques to give a complete answer to the question.\\

	{\bf Acknowledgement}. We would like to thank anonymous referees for their valuable comments and suggestions.

 	\end{document}